\documentclass[xypic,amssymb,rawfonts,amsthm,amsmath,amsfonts,syntonly,11pt]{amsart} 
\usepackage{amsfonts}
\usepackage{amsmath}
\usepackage{amssymb,latexsym}
\usepackage{amsxtra}
\usepackage[matrix,curve,arrow]{xy}

\newtheorem{theorem}{Theorem}[section]

\newtheorem{definition}[theorem]{Definition}

\newtheorem{lemma}[theorem]{Lemma}

\newtheorem{remark}[theorem]{Remark}

\newcommand{\n}{\newline}

\newcommand{\mb}{\mathbb}
\newcommand{\wt}{\widetilde}

\newcommand{\dit}{\raisebox{1ex}{\hspace{0.4mm}.\hspace{0.4mm}}}

\begin{document}

\title[Lefschetz characters for sporadic geometries]
{On non-projective block components of Lefschetz characters for
sporadic geometries}
\date{\today}
\author[R.~Philip Grizzard]{R.~Philip Grizzard}
\thanks{The author was partially supported
        by NSA grant H 98230-05-1-0075.} 

\address{Department of Mathematics (m/c 249),
         University of Illinois at Chicago,
         851 S. Morgan,
         Chicago, IL 60607, USA}
\email{grizzard@math.uic.edu}

\begin{abstract}
This work examines the possible projectivity of 2-modular block
parts of non-projective Lefschetz characters over 2-local
geometries of several sporadic groups.  Previously known results
on $M_{12}$, $J_2$, and $HS$ are mentioned for completeness.  The
main new results are on the sporadic groups $Suz$, $Co_3$, $Ru$,
$O'\!N$, and $H\!e$. For each group, the Lefschetz character is
calculated, and its 2-modular block parts are examined for
projectivity. In each case it is confirmed that a non-principal
block part contains a non-projective summand.  The case of $O'\!N$
is additionally found to have a non-projective summand in its
principal block part. Nineteen of the sporadic groups (including
many previously known cases) are categorized into three classes
based on projectivity properties of their Lefschetz characters.
\end{abstract}

\maketitle

\section{Introduction}\label{intro}

Of the 26 sporadic groups, 11 have 2-local geometries with mod 2
projective Lefschetz modules, and most of these have been
calculated previously \cite{rsy90}. Of the remaining 15 with
non-projective Lefschetz modules, the Lefschetz characters in
terms of irreducibles already have been calculated for three of
these groups. On the other hand, there are 7 of the remaining 12
sporadics whose 2-modular irreducibles are not yet known. Thus
there remain 5 sporadic groups with non-projective Lefschetz
modules whose 2-modular irreducibles are known \cite{dmw}, but
whose Lefschetz characters in characteristic 2 have not yet been
discussed in the literature.

In this work, we calculate the Lefschetz characters and their
2-modular block decompositions for these five sporadic groups,
namely $Suz$, $O'\!N$, $H\!e$, $Co_3$, and $Ru$. We find that the
Lefschetz character for $O'\!N$ has a {\em non}-projective
constituent in its principal block part, while the principal block
parts for the other four groups each has the same character as a
projective module.

This paper is a version of the author's Ph.D. thesis at the
University of Illinois at Chicago.

\section{Notation and Assumed Results}\label{background}

\subsection*{Definitions}\label{lefcharcomp}

The (reduced) {\em Lefschetz module} of a complex $\Delta$ with
action by a group $G$ is the virtual $\mb{Z}G$-module given by the
alternating sum of the chain groups (cf. \cite[p. 281]{rsy90}).
Its degree term is the (reduced) {\em Euler characteristic}:
$\wt{\chi}(\Delta):=\sum^{\dim \Delta}_{i=-1}(-1)^{i}\dim C_i
(\Delta).$  We may reduce this module mod 2, and then analyze it
using the 2-modular representation theory of $G$.

\begin{definition}[Lefschetz character]\label{d:lefchar}
The {\em Lefschetz character} of the natural 2-local
geometry\footnote{``2-local geometry'' is defined in Benson-Smith
\cite[\S 7.2]{bs06}.} $\Delta$ of $G$ is the sequence of Euler
characteristic values of the subgeometries fixed by
representatives $g$ of each conjugacy class of $G$:
$Tr(g,\wt{L})=\wt{\chi}(\Delta^g)$. We denote the Lefschetz
character of $G$ by $\widetilde{\Lambda}_{G}$, or just
$\widetilde{\Lambda}$ when the group is unambiguous.

Let $B$ be a particular 2-modular block for ${G}$. Then we define
$\widetilde{\Lambda}^{G}_{B}$ to be the component of
$\widetilde{\Lambda}_{G}$ coming from Brauer characters in $B$.
\end{definition}

\begin{definition}[Projective]
We define $\Phi(\varphi_i)$ as the {\em character} of the
projective cover $P(\varphi_i)$ of the 2-modular irreducible
Brauer character $\varphi_i$. We say a module has {\em projective
character} if its character can be written as a sum of
$\Phi(\varphi_i)$.
\end{definition}

\begin{definition}[Virtual projective]
\label{d:virtprojmodule} A {\em virtual projective module} is a
$\mb Z$-linear combination of projective indecomposables
$P(\varphi)$.  We say a virtual module has {\em v-projective
character} if the module has the same character as that of a
virtual projective module, i.e. if its character equals some $\mb
Z$-linear combination of $\Phi(\varphi_i)$.  When describing the
character itself with this property, we say the character is {\em
v-projective}.  We emphasize that if a module has v-projective
character, that does not necessarily imply that the module is a
virtual projective module.
\end{definition}

\subsection*{Tests for projectivity}

We use two standard properties of projective modules below.  For
each of these, the contrapositive presents a quick check for
failure of a Lefschetz module to be projective.  The second
condition is stronger.

\begin{lemma}[The $p$-test]\label{lem:ptest}
If $\wt{L}(\Delta)$ is projective, then $|G|_p$ divides
$\wt{\chi}(\Delta)$.
\end{lemma}

\begin{lemma}[The vanishing test]\label{lem:vanishing} \cite[p.148]{feit82}
For the module $\wt{L}(\Delta)$, its Lefschetz character
$\wt{\Lambda}$ is v-projective if and only if $\wt{\Lambda}$
vanishes at the nontrivial 2-elements.
\end{lemma}

If a Lefschetz character passes the vanishing test, this indicates
that the character is v-projective, i.e. can be expressed as a
combination of $\Phi(\varphi_i)$.

\section{Some remarks on the results of our study}\label{remarks}

\subsection*{Partition of the sporadics based on their Lefschetz
characters}\label{sporadicclass}

We define Class I groups as having Lefschetz modules that are
projective. We will refer to each group in this class as a
``Lefschetz Module Projective.''

A Class II group is one whose Lefschetz module is non-projective,
yet whose Lefschetz module restricted to the principal block part
has v-projective character (perhaps zero). This is the case with
three groups previously studied in the literature,
and we will show this to be the case with four more sporadics.
We will refer to a sporadic group in this class as a ``Principal
Block Part V-Projective.''

Finally, a Class III group contains a non-projective summand in
the principle block part of its Lefschetz module---and this part
has no cohomology\footnote{In particular, the variety of this
principal block part must lie in the {\em representation theoretic
nucleus}, as introduced by Benson-Carlson-Robinson \cite[p.
68]{bcr90} (see also Benson \cite[prop. 1.1]{benson02}. } since
Benson-Smith established that the Lefschetz module is acyclic.
Though there are standard constructions in the literature
\cite{cr94} of such modules, this work concretely exhibits (for
the first time) an example of such a module naturally occurring in
the context of $\wt{L}(\Delta)$: in the O'Nan group. We will refer
to a group in this class as a ``Principal Block Part
Non-Projective.''

\subsection*{Other observations}\label{obs}

In section 8.2, Benson and Smith \cite{bs06} remark that reduced
Lefschetz module appears always to involve an indecomposable in a
non-principal block.  After examining larger groups, we can
reformulate their comment as follows.

\begin{remark}[Non-principal block observation]\label{nonprobs}
\label{rem:nonprinc} {\em Each sporadic group affording a
non-projective Lefschetz module appears to have a non-projective
part in a non-principal block.}
\end{remark}

This observation holds in each case we study, and we conjecture
that it will hold for the remaining seven sporadics.

In addition, Smith \cite{sstalk} noticed a further pattern: the
defect of the largest non-principal block part seems to equal the
2-power difference between $|\wt{\chi}(\Delta)|_2$ and $|G|_2$. Of
the sporadic groups we study, only $O'\!N$ fails to follow this
further pattern.\footnote{For an approach to this pattern via
vertices of modules, see Sawabe \cite{sawabe06}.}

\begin{remark}\label{2x}
{\em In these 15 sporadic cases where $\wt{L}(\Delta)$ is not
projective, the vanishing test \ref{lem:vanishing} guarantees that
there will be at least one 2-element $b$ for which $\Delta^b$ is
not contractible. We mention an experimental observation: this
nonzero value for the Euler characteristic often matches the Euler
characteristic value for a well-known non-contractible geometry.
This hints that perhaps our geometric structure $\Delta^b$ is
homotopy equivalent to that well-known geometry. We will suggest
such possible geometries for seven groups in our
study.\footnote{After a preprint version of this paper was
circulated, Silvia Onofrei and John Maginnis \cite{mo07}
established many such homotopy equivalences.} }
\end{remark}

\section{Review of known cases}\label{knowncases}

We present the Lefschetz characters of $M_{12}$, $J_2$, and $HS$,
and then proceed to display their 2-modular block decompositions.

\subsection*{The Mathieu group $M_{12}$}

Benson-Wilkerson \cite[p. 44]{bw95} computed
$\widetilde{\Lambda}_{M_{12}} = \chi_{14}+2\chi_{15}$.  These
characters lie in the unique non-principal block. We express this
as a vector with columns indexed by $M_{12}$ conjugacy classes:
$\widetilde{\Lambda}_{M_{12}}$ = [ $\overbrace{496}^{2^4\cdot 31},
\overset{\underline{\mathbf{2A}}}{\mathbf{-4}},
\overset{\underline{\mathbf{2B}}}{\mathbf{0}}$, -8, -5, 0, 0, 1,
-1, 0, 0, 0, 1, 1, 1 ].

The vanishing test indicates that $\Delta^{2A}$ is not
contractible. From the viewpoint of Remark \ref{2x}, the value -4
in $\widetilde{\Lambda}_{M_{12}}$ matches the degree of the
Lefschetz module of the 
building $\Sigma$ of a subgroup $L_2(4)$ in the centralizer of a
2A element.  This suggests that $\Delta^{2A}$ may be homotopy
equivalent to this $\Sigma$. Finally, one concise way of
expressing the Lefschetz character via possible v-projective and
non-projective parts is: $\widetilde{\Lambda}_{M_{12}} =
\Phi(\varphi_6)+\chi_{15}.$

\subsection*{The Janko group $J_2$}

Benson and Smith \cite[after Theorem 8.7.1]{bs06} found
$\widetilde{\Lambda}_{J_2} = \chi_{16} + \chi_{17} + \chi_{19}$,
noting that these characters lie in the unique non-principal
block. We calculate $\widetilde{\Lambda}_{J_2}$ = [
$\overbrace{736}^{2^5\cdot 23},
\overset{\underline{\mathbf{2A}}}{\mathbf{0}},
\overset{\underline{\mathbf{2B}}}{\mathbf{-4}}$, 16, -5, 0, 1, 1,
-4, -4, 0, -1, 1, 0, 1, 1, 0, 0, 0, 1, 1 ].
Again from the viewpoint of Remark \ref{2x}, the value -4 at class
2B matches the degree of $\wt{L}(\Sigma)$ for building $\Sigma$ of
an $L_2(4)$ in C(2B), the centralizer of conjugacy class 2B. We
can express the Lefschetz character as: $\widetilde{\Lambda}_{J_2}
= \Phi(\varphi_8) + \chi_{16}$.

\subsection*{The Higman-Sims group $HS$}

Klaus Lux calculated the Lefschetz character for $HS$ in a private
communication to Stephen D. Smith: \n $\widetilde{\Lambda}_{HS} =
[\overbrace{42624}^{2^7\cdot 333},
\overset{\underline{\mathbf{2A}}}{\mathbf{0}},
\overset{\underline{\mathbf{2B}}}{\mathbf{-16}}$, -36, 0, 0, 0,
24, -1, 4, 2, 0, 1, 0, 0, 0, 0, -1, -1, -1, 0, -1, 0, 0 ]. In the
viewpoint of Remark \ref{2x}, the value of -16 is the degree of
$\wt{L}(\Sigma)$ for $\Sigma$ the building of $Sp_4(2)$ in C(2B).

\subsubsection*{Block decomposition}
In this case we have some contribution from both blocks. For the
principal block, we find that $\widetilde{\Lambda}^{HS}_{pr}=
2\Phi(\varphi_8),$ so this is v-projective. The ``closest'' we
could come to expressing the block 2 part as a v-projective
character is
$\widetilde{\Lambda}^{HS}_{b2}=\Phi(\varphi_6)+\Phi(\varphi_7)+\Phi(\varphi_9)
+\chi_{24}.$

We separate the Lefschetz character values according to the
2-modular blocks: \noindent $\widetilde{\Lambda}^{HS}_{pr}$ = [
$\overbrace{26624}^{2^{11}\cdot 13},
\overset{\underline{\mathbf{2A}}}{\mathbf{0}},
\overset{\underline{\mathbf{2B}}}{\mathbf{0}}$, -16, 0, 0, 0, 24,
24, 4, 0, 0, -4, 0, 0, 0, 0, 0, 4, 4, 0, -6, 0, 0 ] and
$\widetilde{\Lambda}^{HS}_{b2}$= [ $\overbrace{16000}^{2^7\cdot
125}, \overset{\underline{\mathbf{2A}}}{\mathbf{0}},
\overset{\underline{\mathbf{2B}}}{\mathbf{-16}}$, -20, 0, 0, 0, 0,
-25, 0, 2, 0, 5, 0, 0, 0, 0, -1, -5, -5, 0, 5, 0, 0 ].  We see by
the vanishing test that the non-projective part of
$\widetilde{\Lambda}_{HS}$ is contained completely in the
non-principal block, which is consistent with Remark
\ref{nonprobs}.

\section{The Suzuki group}

We study the maximal subgroups of $Suz$ which contain a Sylow
2-group \cite[p. 131]{atlas}:

\hspace{1in}$H_1\cong 2_-^{1+6}$$\spdot$$U_4(2)$,

\hspace{1in}$H_2\cong 2^{2+8}:(A_5\times S_3)$, and

\hspace{1in}$H_4\cong 2^{4+6}:3A_6$.

\subsection*{Geometry}\label{geometry} The diagram for the natural
geometry of $Suz$ is given in Ronan-Smith \cite[p. 288]{rs80} as:
\[
\begin{xy}<.7cm,0cm>:
(0,0)="a"*{\bullet}; (1,0)="b"*{\bullet}**\dir2{-};
(2,0)="c"*{\bullet}**\dir2{-};
 "a"*+!D{1};"b"*+!D{2};"c"*+!D{4}
\end{xy}
\]
We proceed as in Benson-Smith \cite[Def. 5.3.2]{bs06}. The set
$I:=\{1,2,4\}$ indexes the ``simplex of subgroups''
$\mathcal{H}_I$ of $Suz$. This simplex is determined by the $H_i$,
with $H_J:=\bigcap_{i\in J} H_i$ for all nonempty subsets
$J\subseteq I$.  For simplicity, we will use the notation $H_{ij}$
(without braces) to represent the subgroup $H_{\{ i,j\}}$, etc.

The Lefschetz character in this geometry is an alternating sum of
the permutation representations on vertices (with stabilizers
$H_1$, $H_2$, and $H_4$), edges (with stabilizers $H_{12}$,
$H_{14}$, $H_{24}$), and 2-simplices with stabilizer $H_{124}$.

\subsection*{Inducing from subgroups}

We calculate Lefschetz characters by inducing trivial
representations through a chain of subgroups: from intersections
of $H_1$, $H_2$, and $H_4$ up to the maximal subgroups themselves,
and then up to the whole group.\footnote{This process is described
and several examples are given in section 4 of \cite{rsy90}.}

We then combine these in an alternating sum. Since $H_1$, $H_2$,
and $H_4$ are maximal subgroups of $Suz$, the permutation
characters on their cosets are readily available in GAP
\cite{gap}.\footnote{I am grateful to Professor Robert A. Wilson
of Queen Mary University of London for showing me the step-by-step
details of this process.} The permutation characters for the {\em
intersections} of the maximal subgroups, however, are not
contained in GAP. Instead, selected residues (i.e. links of
simplices) described by the diagram tell us which characters of
the intersections are needed in the calculation. We will
investigate Res($H_1$), in order to induce the trivial characters
of $H_{12}$, $H_{14}$, and $H_{124}$ up to $H_1$.  Then we use
Res($H_4$) to induce from $H_{24}$ up to $H_4$. GAP can then
induce the resulting characters from $H_1$ and $H_4$ up to $Suz$,
as well as the trivial characters for the maximal subgroups
$H_{1}$, $H_{2}$, $H_{4}$ up to $Suz$.

Using the diagram, we find that inducing the trivial character of
$H_{12}$ up to $H_1$ will be the inflation of the trivial
character of $2^4\!\!:\!A_5$ induced up to $U_4(2)$. The Atlas
gives this character to be $1a+6a+20a$. So we have (by the
transitivity of induction):
$$[1a]_{H_{12}}\!\!\uparrow^{Suz}=([1a]_{H_{12}}\!\!\uparrow^{H_1})
_{H_1}\!\!\uparrow^{Suz}=[1a+6a+20a]_{H_1}\!\!\uparrow^{Suz}.$$

The diagram also tells us that $[1a]_{H_{14}}\!\!\uparrow^{H_1}$
is the inflation of the trivial character of $2^{\cdot}(A_4\times
A_4)_{\cdot} 2$ induced up to $U_4(2)$, which is $1a+20a+24a$, so:
$$[1a]_{H_{14}}\!\!\uparrow^{Suz}=([1a]_{H_{14}}\!\!\uparrow^{H_1})
_{H_1}\!\!\uparrow^{Suz}=[1a+20a+24a]_{H_1}\!\!\uparrow^{Suz}.$$

\noindent Similarly $[1a]_{H_{124}}\!\!\uparrow ^{H_{12}}$ is the
inflation of $[1a]_{A_4}\!\!\uparrow ^{A_5}=1a+4a$, giving us:
$$[1a]_{H_{124}}\!\!\uparrow^{Suz}=(([1a]_{H_{124}}\!\!\uparrow
^{H_{12}})_{H_{12}}\!\!\uparrow^{H_1})_{H_1}\!\!\uparrow^{Suz}=
([1a+4a]_{H_{12}} \!\!\uparrow^{H_1})_{H_1}\!\!\uparrow^{Suz}.$$

There is a mild complication here, in that GAP does not contain
the character tables of maximal subgroups such as $2^4\!\!:\!A_5$
of $U_4(2)$, so it is not set up to do the induction
$[4a]_{H_{12}} \!\!\uparrow^{H_1}$.  But this can be calculated
using the viewpoint of Smith \cite{sscomac}. We display the result
in the style of the Atlas:
$$[1a]_{H_{124}}\!\!\uparrow^{Suz}=([1a]_{H_{124}}\!\!\uparrow^{H_1})
_{H_1}\!\!\uparrow^{Suz}=[1a+6a+20aa+24a+64a]_{H_1}\!\!\uparrow^{Suz}.$$


Res($H_1$) has taken care of all but one of the intersections of
maximal subgroups whose trivial characters we need to induce up to
$Suz$. The last one needed is $H_{24}$, which we examine in
Res($H_2$). We have:
$$[1a]_{H_{24}}\!\!\uparrow^{Suz}=([1a]_{H_{24}}\!\!\uparrow^{H_2})
_{H_2}\!\!\uparrow^{Suz}=[1a+4a]_{H_2}\!\!\uparrow^{Suz}.$$


We are now ready to calculate the Lefschetz character of $Suz$ as
an alternating sum, with the sign at a simplex $\sigma$ given by
$(-1)^{\dim\sigma}$. To summarize our above findings:
\begin{align}
\phantom{-}[1a]_{H_{124}}\!\!\uparrow^{Suz} &= \phantom{-}[1a + 6a
+ 20aa + 24a +64a]_{H_1}\!\!\uparrow^{Suz} \notag \\
-[1a]_{H_{12}}\!\!\uparrow^{Suz} &= -[1a + 6a + 20a]
_{H_1}\!\!\uparrow^{Suz}\notag \\
-[1a]_{H_{14}}\!\!\uparrow^{Suz} &= -[1a + 20a + 24a]
_{H_1}\!\!\uparrow^{Suz}\notag \\
\phantom{-}[1a]_{H_1}\!\!\uparrow^{Suz} &= \phantom{-}
[1a]_{H_1}\!\!\uparrow^{Suz}\notag \\
\hline -[1a]_{H_{24}}\!\!\uparrow^{Suz} &= -[1a + 4a]
_{H_2}\!\!\uparrow^{Suz}\notag \\
\phantom{-}[1a]_{H_2}\!\!\uparrow^{Suz} &= \phantom{-}
[1a]_{H_2}\!\!\uparrow^{Suz}\notag \\
\hline \phantom{-}[1a]_{H_4}\!\!\uparrow^{Suz} &=
\phantom{-}[1a]_{H_4}\!\!\uparrow^{Suz}\notag
\end{align}
We can do some quick cancellation to simplify:
\begin{equation}\label{E:LefSuz}
\widetilde{\Lambda}_{Suz} = [64a]_{H_1}\!\!\uparrow^{Suz} -
[4a]_{H_2}\!\!\uparrow^{Suz} +
[1a]_{H_4}\!\!\uparrow^{Suz}-1.\notag
\end{equation}

\noindent The result is: $\widetilde{\Lambda}_{Suz}$ = [4189184,
0, -64, 3968, -352, -73, 0, 0, 0, 8, -16, 9, 0, 0, 0, 0, -1, -1,
0, 0, 0, 2, 2, 0, 1, -1, 0, 0, 0, -1, 0, -1, -1, -1, 2, 2, 3, 0,
0, 0, -1, -1, 0]. {\em This is the Lefschetz character in vector
form, i.e., the columns are indexed by the conjugacy classes.}

\subsection*{Block decomposition}

As a sum of irreducible complex characters, we get
$\widetilde{\Lambda}_{Suz}$=[ 0, 0, 0, 0, 0, 0, 0, 0, 0, 0, 0, 0,
0, 0, 0, 0, 0, 0, 0, 0, 0, 0, 0, 2, 1, 1, 1, 0, 1, 0, 1, 1, 0, 0,
2, 2, 1, 2, 2, 2, 2, 2, 3 ]. {\em This is the Lefschetz character
in scalar product form, i.e., the columns are indexed by the
irreducible complex characters of $Suz$.}

We use these coefficients of irreducible characters along with the
decomposition matrix to decompose $\widetilde{\Lambda}_{Suz}$ into
v-projective and non-projective parts. We find that
$\widetilde{\Lambda}^{Suz}_{pr}= \Phi(\varphi_{13}) +
\Phi(\varphi_{14})$, a v-projective character.

For block 2 (of defect 3), the ``closest'' we can come to
decomposing $\widetilde{\Lambda}^{Suz}_{b2}$ into v-projectives
is:
$\widetilde{\Lambda}^{Suz}_{b2}=\Phi(\varphi_{15})+\Phi(\varphi_{16})+
\chi_{43}.$ Hence the block 2 part of the Lefschetz character of
$Suz$ is non-projective, verifying our observation in
Remark~\ref{nonprobs}. $Suz$ is categorized into Class II.
Furthermore, the non-principal block part has defect 3, which
matches the $2$-power difference between $|G|_2=2^{13}$ and
$|\widetilde{\chi}(\Delta)|_2=|\widetilde{\Lambda}_{Suz}|_2=2^{10}$,
continuing the pattern noticed by Smith \cite{sstalk}.

We summarize our findings with the following theorem.
\begin{theorem}\label{thm:suz}
$\widetilde{\Lambda}^{Suz}_{pr}= \Phi(\varphi_{13}) +
\Phi(\varphi_{14})$. Thus $Suz$ is a Class II ``Principal Block
Part V-Projective'' group.
\end{theorem}

\subsection*{The Lefschetz character values separated into
blocks}

The vectors of Lefschetz character values (indexed by the
conjugacy classes) given just by characters in each respective
block are:\n $\widetilde{\Lambda}^{Suz}_{pr}$ = [ 2834432,
$\overset{\underline{\mathbf{2A}}}{\mathbf{0}}$,
$\overset{\underline{\mathbf{2B}}}{\mathbf{0}}$, 512, -352, -64,
0, 0, 0, 0, 32, 32, 0, 0, 0, 0, 0, -8, 0, 0, 0, -4, -4, 0, 0, -4,
0, 0, 0, 0, 0, 3, 3, 0, -4, -4, 2, 0, 0, 0, 1, 1, 0 ] and
\noindent $\widetilde{\Lambda}^{Suz}_{b2}$= [ 1354752,
$\mathbf{0}$, $\mathbf{-64}$, 3456, 0, -9, 0, 0, 0, 8, -48, -23,
0, 0, 0, 0, -1, 7, 0, 0, 0, 6, 6, 0, 1, 3, 0, 0, 0, -1, 0, -4, -4,
-1, 6, 6, 1, 0, 0, 0, -2, -2, 0 ]. Here we have used bold face for
the values at elements of order 2 in order to make clear the
results of the vanishing test \ref{lem:vanishing}. From the
viewpoint of Remark~\ref{2x}, this value of -64 is also the degree
of $\wt{L}(\Sigma)$ for the building $\Sigma$ of $L_3(4)$ in
C(2B), suggesting a possibly homotopy equivalence.

\section{The O'Nan group}

For our remaining four groups starting with $O'\!N$, we present a
somewhat more abbreviated study with fewer details than with
$Suz$. The maximal subgroups of interest are \cite[p. 132]{atlas}:

\hspace{1in} $H_1\cong 4_2\dit L_3(4)\!\!:\!2_1$, and

\hspace{1in} $H_3\cong 4^3\dit L_3(2)$.

With only two maximal subgroups, the process of obtaining the
Lefschetz character is much simpler.  There is no need for a
diagram, and we have only one intersection to find: $H_{13}$. We
can induce the trivial character of $H_{13}$ via:
$$[1a]_{H_{13}}\!\!\uparrow^{O'\!N}=([1a]_{H_{13}}\!\!\uparrow^{H_3})
_{H_3}\!\!\uparrow^{O'\!N}=[1a+6a]_{H_3}\!\!\uparrow^{O'\!N}.$$


One easy cancellation leads to the final calculation:
\begin{equation}\label{E:LefON}
\widetilde{\Lambda}_{O'\!N} = [1a]_{H_1}\!\!\uparrow^{O'\!N} -
[6a]_{H_3}\!\!\uparrow^{O'\!N}-1.\notag
\end{equation}

Using GAP, we obtain the Lefschetz character (negating to obtain
positive degree) $\widetilde{\Lambda}_{O'\!N}$ =
[$\overbrace{254294272}^{2^8\cdot \,993337}, \overset{\underline{\
\mathbf{2A}\ }}{\mathbf{8960}}$, -44, 0, 0, -8, -4, -48, -13, 0,
0, 0, 1, 0, 0, 1, 1, 0, 0, 0, 0, 1, 1, 1, 0, 0, 0, 0, 1, 1].

\subsection*{The principal block} We might expect that the
principal block part has v-projective character as it did in
$Suz$.  But we find:

$\widetilde{\Lambda}^{O'\!N}_{pr}= [
\overbrace{147299328}^{2^{10}\cdot 143847}, \overset{\underline{\
\mathbf{2A}\ }}{\mathbf{7168}}$, -2736, 0, 0, 708, -32, 2388, 344,
0, 0, 28, 176, 0, -28, -156, -156, 0, 0, 0, 0, 80, 80, 80, 0, 0,
0, 0, -520, -520 ].

Notice that $|O'\!N|_2=2^9$ does indeed divide
$|\widetilde{\Lambda}^{O'\!N}_{pr}|$, so the $p$-test
\ref{lem:ptest} passes.

However, the vanishing test \ref{lem:vanishing} fails, and so
$\widetilde{\Lambda}^{O'\!N}_{pr}$ is actually {\em not}
projective. Indeed this block part is in some sense very far from
being v-projective, as the ``closest'' expression (in terms of
minimizing the coefficients of non-projective components) we can
find is:
\begin{center}
\parbox{4in}{\begin{align}
\widetilde{\Lambda}^{O'\!N}_{pr} &=
20\Phi(\varphi_3)+12\Phi(\varphi_7)+
16\Phi(\varphi_8) \notag\\
  &+ 2\chi_8+2\chi_9+18\chi_{11}+
4\chi_{12}+4\chi_{13}+4\chi_{14}+2\chi_{15} \notag\\
  &+ 8\chi_{18}+
6\chi_{19}+12\chi_{20}+18\chi_{21}+18\chi_{22}+14\chi_{25}. \notag
\end{align}}
\end{center}
\subsubsection*{Block 2---defect 3} We compute
$\widetilde{\Lambda}^{O'\!N}_{b2}=[ \overbrace{1386240}^{2^8\cdot
5415}, \overset{\underline{\ \mathbf{2A}\ }}{\mathbf{1792}}$,
1140, 0, 0, 60, 28, 324, -12, 0, 0, -28, -24, 0, 28, 60, 60, 0, 0,
0, 0, 0, 0, 0, 0, 0, 0, 0, -18, -18]. This block part is also
non-projective, and we can write
$\widetilde{\Lambda}^{O'\!N}_{b2}= 8\Phi(\varphi_2) +
4\Phi(\varphi_4) + 8\Phi(\varphi_5) + 14\chi_7.$

\subsection*{The defect 0 blocks}

The Lefschetz character of $O'\!N$ has several constituents from
blocks of defect 0, which are of course projective:
$\widetilde{\Lambda}^{O'\!N}_{b3}=97\varphi_9$,
$\widetilde{\Lambda}^{O'\!N}_{b4}=97\varphi_{10}$,
$\widetilde{\Lambda}^{O'\!N}_{b5}=115\varphi_{11}$,
$\widetilde{\Lambda}^{O'\!N}_{b6}=115\varphi_{12}$,
$\widetilde{\Lambda}^{O'\!N}_{b7}=115\varphi_{13}$.  These do not
seem very illuminating for our study.

\begin{theorem}\label{thm:on}
By the vanishing test \ref{lem:vanishing},
$\widetilde{\Lambda}^{O'\!N}_{pr}$ is non-projective. Thus $O'\!N$
is categorized into Class III---it is a ``Principal Block Part
Non-Projective'' sporadic group.
\end{theorem}

\subsection*{Observations of unusual behavior}\label{onobs}

The principal block part of the Lefschetz character of $O'\!N$
contains a non-projective summand, making $O'\!N$ the charter
member of Class III. This discovery for O'Nan provides the first
``natural'' example in the literature on sporadic geometries of
such a principal block part with no cohomology.

The principal block part of $O'\!N$ is the only block part we find
in our study to pass the $p$-test \ref{lem:ptest}, yet fail the
vanishing test \ref{lem:vanishing}.  This may be related to
$O'\!N$ being the only group in our study with only one conjugacy
class of order 2. Furthermore, the value of
$\widetilde{\Lambda}_{O'\!N}$ at that involution is 8960.  This is
an unfamiliar Euler characteristic value, and indeed $O'\!N$ is
the only group for which we do not make a homotopy equivalence
conjecture as in Remark~\ref{2x}.

Finally, note that the largest non-principal block of $O'\!N$ has
defect 3 but the 2-power difference between $|\wt{\Lambda}|_2$ and
$|G|_2$ is only 1.  This is the only example we study that
violates the pattern observed by Smith \cite{sstalk}.

\section{The Held group}

Our next group is the Held group $H\!e$. The relevant maximal
2-local subgroups are \cite[p. 131]{atlas}:

\hspace{1in}$H_1\cong 2_+^{1+6}.L_3(2)$,

\hspace{1in}$H_{6a}\cong 2^{6}:3$$\spdot$$S_6$, and

\hspace{1in}$H_{6b}\cong 2^{6}:3$$\spdot$$S_6$.

The diagram is given in Ronan-Smith \cite[p. 288]{rs80} as:
\[\begin{xy}<.7in,0in>:
(0,0)="a"*{\bullet}; (1,0)="b"*{\bullet}**\dir{-};
(.5,-\halfrootthree)="c"*{\circ}**\dir2{-};"a"**\dir2{-};
"a"*+!D{6a}; "b"*+!D{6b}; "c"*+!U{1};
\end{xy}\]
We obtain the subgroups

\hspace{1cm}$H_{1,6a}=H_{1,6b}=S_32^22_+^{1+6}$,
$H_{1,6a,6b}=\frac{\frac{2}{2^2}}{2_+^{1+6}}$, and
$H_{6a,6b}=\frac{3(S_4\times 2)}{2^6}$.

Using previously described techniques, we find:\footnote{I am
grateful to Silvia Onofrei and John Maginnis for correcting an
error at this point in my original calculation.}
\begin{equation}\label{E:LefHe}
\widetilde{\Lambda}_{H\!e} = [8a]_{H_1}\!\!\uparrow^{H\!e} - [5d +
9b]_{H_{6a}}\!\!\uparrow^{H\!e}+
[1a]_{H_{6b}}\!\!\uparrow^{H\!e}-1.\notag
\end{equation}

With columns indexed by the conjugacy classes of $H\!e$, we
calculate: \n $\widetilde{\Lambda}_{H\!e}$ = [
$\overbrace{1120384}^{2^7\cdot 8753},
\overset{\underline{\mathbf{2A}}}{\mathbf{-64}},
\overset{\underline{\mathbf{2B}}}{\mathbf{0}}$, -197, -8, -8, 0,
0, 9, -1, 0, -1, -1, -1, 6, 6, 0, 1, 1, 0, -1, -1, 0, 0, 3, -1,
-1, -1, -1, -1, -1, -1, -1 ].

\subsection*{Block decompositions}

\vspace{.25cm}$\widetilde{\Lambda}^{H\!e}_{pr}=[
\overbrace{760832}^{2^{10}\cdot 743},
\overset{\underline{\mathbf{2A}}}{\mathbf{0}},
\overset{\underline{\mathbf{2B}}}{\mathbf{0}}$, 152, -40, 0, 0, 0,
-68, 0, 0, 16, 16, 58, -12, -12, 0, 0, 0, 0, 0, 0, 0, 0, -8, 14,
14, -2, -2, -5, -5, 0, 0 ]. Using the 2-modular table, we can
express this as
$\widetilde{\Lambda}^{H\!e}_{pr}=\Phi(\varphi_{10})+
\Phi(\varphi_{11})+ \Phi(\varphi_{12})$. For block 2 (of defect
3), we have \vspace{.25cm}$\widetilde{\Lambda}^{H\!e}_{b2}= [
\overbrace{58496}^{2^7\cdot 457},
\overset{\underline{\mathbf{2A}}}{\mathbf{-64}},
\overset{\underline{\mathbf{2B}}}{\mathbf{0}}$, -13, 32, -8, 0, 0,
21, -1, 0, -17, -17, 39, 18, 18, 0, 1, 1, 0, -1, -1, 0, 0, -3, -1,
-1, -6, -6, 4, 4, -1, -1 ]. This fails tests \ref{lem:ptest} and
\ref{lem:vanishing}, so we use the decomposition matrix to see how
``close'' this block part is to being v-projective:
$\widetilde{\Lambda}^{H\!e}_{b2}= 3\Phi(\varphi_{14}) +
\chi_{15}$. From the viewpoint of Remark \ref{2x}, this value of
-64 is also the degree of $\wt{L}(\Sigma)$ for the building
$\Sigma$ of $L_3(4)$ in C(2B). The parts in blocks 3 and 4 must be
projective, since they have defect 0. We have
$\widetilde{\Lambda}^{H\!e}_{b3}=7\varphi_{15}$, and
$\widetilde{\Lambda}^{H\!e}_{b4}=7\varphi_{16}$.

\begin{theorem}\label{thm:he}
$\widetilde{\Lambda}^{H\!e}_{pr}=\Phi(\varphi_{10})+
\Phi(\varphi_{11})+ \Phi(\varphi_{12})$.  Hence $H\!e$ is a Class
II ``Principal Block Part V-Projective'' group.
\end{theorem}

\section{The Conway group $Co_3$}

The maximal subgroups we examine are \cite[p. 134]{atlas}:

\hspace{1in}$H_1\cong 2\spdot Sp_6(2)$,

\hspace{1in}$H_2\cong 2^2.\,[2^7.\,3^2].\,S_3$, and

\hspace{1in}$H_4\cong 2^4$$\spdot$$A_8$.

The diagram for the natural geometry \cite[\S8.13]{bs06} is:
\[\begin{xy}<.7cm,0cm>:
(0,0)="a"*{\bullet}; (1,0)="b"*{\bullet}**\dir2{-};
(2,0)="c"*{\bullet}**\dir1{-};"b";(1,-1)*+{\Box}**\dir{-};
 "a"*+!D{1};"b"*+!D{2};"c"*+!D{4}
\end{xy}\]

This gives us the needed intersections
$H_{12}=2^2.[2^6]\!\!:\!(S_3\times S_3)$,
$H_{14}=2.2^6\!\!:\!L_3(2)$, $H_{124}=2.2^6\!\!:\!S_4$, and
$H_{24}=2^4.2^4\!\!:\!(S_3\times S_3)$.  These lead us to:
\begin{equation}\label{E:LefCo3}
\widetilde{\Lambda}_{Co_3} = [216a+280b]_{H_1}\!\!\uparrow^{Co_3}
- [14a+20a]_{H_4}\!\!\uparrow^{Co_3} +
[1a]_{H_2}\!\!\uparrow^{Co_3}-1.\notag
\end{equation}

Now we use GAP to obtain $\widetilde{\Lambda}_{Co_3}$ = [
$\overbrace{50378624}^{2^7\cdot 393583},
\overset{\underline{\mathbf{2A}}}{\mathbf{0}},
\overset{\underline{\ \mathbf{2B}\ }}{\mathbf{-496}}$, -2080,
-784, 125, 0, 0, 24, 19, 0, 0, 0, 8, 5, 2, 0, 0, 0, 8, -1, 0, -1,
-1, -1, 0, 0, 0, 0, 0, 1, 0, 0, 0, -1, -1, -1, -1, -1, 0, 0, 0 ].
From the viewpoint of Remark \ref{2x}, we note that -496 is the
degree of $\wt{L}(\Sigma)$ for the 2-local geometry $\Sigma$ of an
$M_{12}$ in C(2B).

\subsection*{The Principal block} We compute
$\widetilde{\Lambda}^{Co_3}_{pr}$ = [
$\overbrace{34263040}^{2^{12}\cdot 8365},
\overset{\underline{\mathbf{2A}}}{\mathbf{0}},
\overset{\underline{\mathbf{2B}}}{\mathbf{0}}$, -11840, 544, -512,
0, 0, 440, -40, 0, 0, 0, 0, 0, -56, 0, 0, 0, 28, 106, 0, 0, 20,
20, 0, 0, 0, 0, -40, 44, 0, 0, 0, -8, 0, 0, -37, -37, 0, 0, 0 ],
which is v-projective:
$\widetilde{\Lambda}^{Co_3}_{pr}=\Phi(\varphi_9)+\Phi(\varphi_{10})+
6\Phi(\varphi_{12})+ 8\Phi(\varphi_{14}).$

\subsection*{Block 2---defect 3} $\widetilde{\Lambda}^{Co_3}_{b2}$
 = [ $\overbrace{13006720}^{2^7\cdot 101615},
\overset{\underline{\mathbf{2A}}}{\mathbf{0}},
\overset{\underline{\ \mathbf{2B}\ }}{\mathbf{-496}}$, 11296,
-2384, 445, 0, 0, -80, 155, 0, 0, 0, 8, 5, 34, 0, 0, 0, 4, -83, 0,
-1, -21, -21, 0, 0, 0, 0, 16, -19, 0, 0, 0, -17, -1, -1, 36, 36,
0, 0, 0 ]. This character fails either test. We can express this
as: $\widetilde{\Lambda}^{Co_3}_{b2}=
4\Phi(\varphi_{11})+6\Phi(\varphi_{13})+
11\Phi(\varphi_{16})+\chi_{32}+2\chi_{38}.$

\subsection*{Block 3---defect 1} $\widetilde{\Lambda}^{Co_3}_{b3}$
 = [ $\overbrace{3108864}^{2^{12}\cdot 759},
\overset{\underline{\mathbf{2A}}}{\mathbf{0}},
\overset{\underline{\mathbf{2B}}}{\mathbf{0}}$, -1536, 1056, 192,
0, 0, -336, -96, 0, 0, 0, 0, 0, 24, 0, 0, 0, -24, -24, 0, 0, 0, 0,
0, 0, 0, 0, 24, -24, 0, 0, 0, 24, 0, 0, 0, 0, 0, 0, 0 ]. We can
express this as
$\widetilde{\Lambda}^{Co_3}_{b3}=12\Phi(\varphi_{15}).$

\begin{theorem}\label{thm:co3}
$\widetilde{\Lambda}^{Co_3}_{pr}=\Phi(\varphi_9)+\Phi(\varphi_{10})+
6\Phi(\varphi_{12})+ 8\Phi(\varphi_{14})$, making $Co_3$ a Class
II ``Principal Block Part V-Projective'' group.
\end{theorem}

\section{The Rudvalis group}

The relevant maximal 2-local subgroups are \cite[p. 126]{atlas}:

\hspace{1in}$H_1\cong 2\spdot 2^{4+6}\!\!:\!S_5$,

\hspace{1in}$H_3\cong 2^{3+8}\!\!:\!L_3(2)$, and

\hspace{1in}$H_6\cong (2^6\!\!:\!U_3(3))\!\!:\!2\cong 2^6 G_2(2)$.
The diagram is:
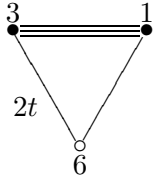
\[\begin{xy}<.7in,0in>:
(0,0)="a"*{\bullet}; (1,0)="b"*{\bullet}**\dir3{-};
(.5,-\halfrootthree)="c"*{\circ}**\dir{-};"a"**\dir{-};
?*+!UR{2t}; "a"*+!D{3}; "b"*+!D{1}; "c"*+!U{6};
\end{xy}\]
in the notation of \cite[p. 288]{rs80}. The diagram gives us the
needed intersections, of the following structures:

\begin{center} $H_{16}=\frac{\frac{S_3}{2^{2+1+2}}}{2^6}$,
$H_{36}=\frac{\frac{S_3}{2^{1+4}}}{2^6}$,
$H_{136}=\frac{\frac{2}{2^{1+4}}}{2^6}$, and
$H_{13}=\frac{\frac{S_3}{2^2}}{2^{3+8}}$.\end{center}

These subgroups allow us to calculate:
\begin{equation}\label{E:LefRu}
\widetilde{\Lambda}_{Ru} = [32ab]_{H_6}\!\!\uparrow^{Ru} -
[6a]_{H_3}\!\!\uparrow^{Ru} + [1a]_{H_1}\!\!\uparrow^{Ru}-1.\notag
\end{equation}

\noindent Indexed by conjugacy classes, we obtain:
$\widetilde{\Lambda}_{Ru}$= [ $\overbrace{10113024}^{2^{12}\cdot
2469}, \overset{\underline{\mathbf{2A}}}{\mathbf{0}},
\overset{\underline{\mathbf{2B}}}{\mathbf{64}}$, -96, 0, 0, 0, 0,
24, -1, 0, 5, 0, 0, 0, 0, -1, 0, 0, -1, 1, 1, 1, -1, 0, 0, 0, 0,
0, 0, 0, -1, -1, -1, -1, -1 ]. This character does not vanish at
the 2B element, telling us that the Lefschetz module is
non-projective.  From the perspective of Remark \ref{2x}, we see
that 64 matches the degree of $\wt{L}(\Sigma)$, where $\Sigma$ is
the building of $Sz(8)$ in C(2B).

\subsection*{Block decomposition}$\widetilde{\Lambda}^{Ru}_{pr}$ = [
$\overbrace{6881280}^{2^{16}\cdot 105},
\overset{\underline{\mathbf{2A}}}{\mathbf{0}},
\overset{\underline{\mathbf{2B}}}{\mathbf{0}}$, -48, 0, 0, 0, 0,
280, 80, 0, 28, 0, 0, 0, 0, 0, 0, 0, -36, 0, 0, 0, 2, 0, 0, 0, 0,
0, 0, 0, 0, 0, 0, -14, -14 ], which is v-projective:
$\widetilde{\Lambda}^{Ru}_{pr}=\Phi(\varphi_5)+\Phi(\varphi_8)$.
Thus $Ru$ belongs in Class II.

For block 2 (of defect 2), we find $\widetilde{\Lambda}^{Ru}_{b2}$
= [ $\overbrace{3231744}^{2^{12}\cdot 789},
\overset{\underline{\mathbf{2A}}}{\mathbf{0}},
\overset{\underline{\mathbf{2B}}}{\mathbf{64}}$, -48, 0, 0, 0, 0,
-256, -81, 0, -23, 0, 0, 0, 0, -1, 0, 0, 35, 1, 1, 1, -3, 0, 0, 0,
0, 0, 0, 0, -1, -1, -1, 13, 13 ]. This does not vanish at element
2B, so $\widetilde{\Lambda}^{Ru}_{b2}$ is non-projective. We
minimize the number of non-projective characters and express the
Lefschetz character of this block part as
$\widetilde{\Lambda}^{Ru}_{b2}=\Phi(\varphi_6)+\Phi(\varphi_7)+
6\Phi(\varphi_9)+\chi_{36}.$

\begin{theorem}\label{thm:ru}
$\widetilde{\Lambda}^{Ru}_{pr}=\Phi(\varphi_5)+\Phi(\varphi_8)$.
So $Ru$ is a Class II ``Principal Block Part V-Projective'' group.
\end{theorem}

\section{Summary}

We present a partition of the sporadic groups according to our
findings, as described in Section \ref{sporadicclass} (boldface
indicates groups originally classified by this work).

\begin{table}[h]
\begin{center}
\begin{tabular}{|p{6.65cm}|p{5.4cm}|}
\hline \textbf{Class} & \textbf{Sporadic Groups}\\
\hline I: Lefschetz Module Projective & $M_{11}$, $J_1$, $M_{22}$,
$M_{23}$, $J_3$, $M_{24}$, $M^cL$, $Co_2$, $Ly$, $J_4$, $Th$ \\
\hline II: Principal Block Part V-Projective & $M_{12}$, $J_2$,
$HS$, $\mathbf{Suz}$, $\mathbf{H\!e}$, $\mathbf{Co_3}$,
$\mathbf{Ru}$ \\
\hline III: Principal Block Part Non-Projective &
$\mathbf{O'\!N}$\\
\hline
\end{tabular}
\label{tab:sporadicclass}
\end{center}
\end{table}

The groups $Co_1$, $Fi_{22}$, $Fi_{23}$, $Fi_{24}'$, $HN$, $B$,
and $M$ have yet to be classified.\footnote{Maginnis and Onofrei
\cite{mo07} have recently classified $Fi_{22}$ into Class II.} The
2-modular decomposition matrices \cite{dmw} of these groups are
not yet known. We do know by the $p$-test that their Lefschetz
modules are non-projective, so they will be in Class II or Class
III.

\section{Future directions}

For the block parts of Lefschetz modules that we found to have
v-projective character, we suggest that the method of Steinberg
module inversion by Webb \cite{webb87b} as extended by Grodal
\cite{grodal02} could possibly be used to show that the block part
is actually a projective module in each case. The case of $O'\!N$
should be investigated further; we found a number of unusual
features (see Section \ref{onobs}) that made this group stand out.
There are still seven sporadic groups whose Lefschetz characters
could be computed (but they could not be fully decomposed mod 2
since the 2-modular irreducibles of these groups are not yet
known). For even further directions, some different geometries for
each group could be studied, as indicated by Benson-Smith
\cite{bs06}. The block decomposition corresponding to primes
$p\neq2$ could also be explored.

\bibliographystyle{amsalpha}

\bibliography{Thesis}

\newcommand{\etalchar}[1]{$^{#1}$}
\providecommand{\bysame}{\leavevmode\hbox to3em{\hrulefill}\thinspace}
\providecommand{\MR}{\relax\ifhmode\unskip\space\fi MR }
\providecommand{\MRhref}[2]{%
  \href{http://www.ams.org/mathscinet-getitem?mr=#1}{#2}
}
\providecommand{\href}[2]{#2}
\begin{thebibliography}{CCN{\etalchar{+}}85}

\bibitem[BCR90]{bcr90}
D.~J. Benson, J.~F. Carlson, and G.~R. Robinson, \emph{On the vanishing of
  group cohomology}, J. Algebra \textbf{131} (1990), no.~1, 40--73.
  \MR{MR1054998 (91c:20073)}

\bibitem[Ben02]{benson02}
D.~J. Benson, \emph{The nucleus and extensions between modules for a finite
  group}, Representations of algebra. Vol. I, II, Beijing Norm. Univ. Press,
  Beijing, 2002, Easiest to find at
  \verb!http://www.maths.abdn.ac.uk/~bensondj/html/archive/benson.html!,
  pp.~145--155. \MR{MR2067376 (2005c:20018)}

\bibitem[Bre99]{dmw}
Thomas Breuer, \emph{Decomposition {M}atrices}, See details at
  \verb!http://www.math.rwth-aachen.de/~MOC/decomposition/!, 1999.

\bibitem[BS07]{bs06}
David~J. Benson and Stephen~D. Smith, \emph{Classifying spaces of sporadic
  groups}, Surveys and Monographs of the AMS, American Mathematical Society,
  2007, Submitted for publication. Preprint available at
  \verb!http://www.maths.abdn.ac.uk/~bensondj/html/archive/benson-smith.html!

\bibitem[BW95]{bw95}
D.~J. Benson and C.~W. Wilkerson, \emph{Finite simple groups and {D}ickson
  invariants}, Homotopy theory and its applications (Cocoyoc, 1993), Contemp.
  Math., vol. 188, Amer. Math. Soc., Providence, RI, 1995, pp.~39--50.
  \MR{MR1349127 (96d:55009)}

\bibitem[CCN{\etalchar{+}}85]{atlas}
J.H. Conway, R.T. Curtis, S.P. Norton, R.A. Parker, and R.A. Wilson,
  \emph{Atlas of finite groups}, Oxford University Press, 1985.

\bibitem[CR94]{cr94}
Jon~F. Carlson and Geoffrey~R. Robinson, \emph{Varieties and modules with
  vanishing cohomology}, Math. Proc. Cambridge Philos. Soc. \textbf{116}
  (1994), no.~2, 245--251. \MR{MR1281544 (95c:20073)}

\bibitem[Fei82]{feit82}
Walter Feit, \emph{The representation theory of finite groups}, North-Holland
  Mathematical Library, vol.~25, North-Holland Publishing Co., Amsterdam, 1982.
  \MR{MR661045 (83g:20001)}

\bibitem[Gro99]{gap}
The~{GAP} Group, \emph{Gap --- {G}roups, {A}lgorithms, and {P}rogramming,
  {V}ersion 4.2}, (\verb!http://www-gap.dcs.st-and.ac.uk/~gap!) (Aachen, St.
  Andrews), 1999.

\bibitem[Gro02]{grodal02}
Jesper Grodal, \emph{Higher limits via subgroup complexes}, Ann. of Math. (2)
  \textbf{155} (2002), no.~2, 405--457. \MR{MR1906592 (2003g:55025)}

\bibitem[MO07]{mo07}
John Maginnis and Silvia Onofrei, \emph{Fixed point sets of involutions},
  unpublished (2007), Preprint available at
  \verb!http://www.math.ksu.edu/~onofrei/fixedpoint.pdf!

\bibitem[RS80]{rs80}
M.A. Ronan and S.~D. Smith, \emph{2-local geometries for some sporadic groups},
  Proceedings of Symposia in Pure Mathematics \textbf{37} (1980), 283--289.

\bibitem[RSY90]{rsy90}
A.~J.~E. Ryba, Stephen~D. Smith, and Satoshi Yoshiara, \emph{Some projective
  modules determined by sporadic geometries}, Journal of Algebra \textbf{129}
  (1990), 279--311.

\bibitem[Saw06]{sawabe06}
Masato Sawabe, \emph{On the reduced {L}efschetz module and the centric
  {$p$}-radical subgroups. {II}}, J. London Math. Soc. (2) \textbf{73} (2006),
  no.~1, 126--140. \MR{MR2197374 (2006j:20078)}

\bibitem[Smi90]{sscomac}
Stephen~D. Smith, \emph{On decomposition of modular representations from
  {C}ohen-{M}acaulay geometries}, Journal of Algebra \textbf{131} (1990),
  598--625.

\bibitem[Smi05]{sstalk}
\bysame, \emph{Cohomology decomposition from subgroup complexes of finite
  groups}, See details at \verb!http://www.math.uic.edu/~smiths/talk.pdf!,
  2005.

\bibitem[Web87]{webb87b}
P.~J. Webb, \emph{A local method in group cohomology}, Comment. Math. Helv.
  \textbf{62} (1987), no.~1, 135--167. \MR{MR882969 (88h:20065)}

\end{thebibliography}

\end{document}